\documentclass[11pt]{article}
\usepackage[english]{babel}
\usepackage{graphics}
\usepackage{amsfonts}
\usepackage{amsmath}
\usepackage{amssymb}
\usepackage{mathrsfs}
\usepackage{enumerate}
\usepackage{epsf}
\usepackage{epsfig}
\usepackage{amscd}
\usepackage{color}
\usepackage[unicode]{hyperref}

\hypersetup{
    colorlinks = true,%
    citecolor = [rgb]{0.0,0.0,0.9},
    filecolor=black,%
    linkcolor = [rgb]{0.65,0.0,0.0},%
    anchorcolor = red,
    pagecolor = red,
    urlcolor= [rgb]{0.65,0.0,0.0}
}

\newtheorem{lemma}{Lemma}[section]

\newtheorem{proposition}{Proposition}[section]
\newtheorem{example}{Example}[section]
\newtheorem{theorem}{Theorem}[section]

\newtheorem{remark}{Remark}[section]
\newtheorem{definition}{Definition}[section]
\def\psn{\par\smallskip\noindent}
\def\QED{\hfill $\Box$\par\smallskip\noindent}

\def\scatola{\lower5pt\hbox{\vbox{\hrule\hbox{\vrule\kern2pt\vbox%
{\kern5pt\hbox{\mathsurround=0pt }\kern2pt}\kern4pt\vrule}\hrule}}\
} %%
\def\l{\lambda}

\def\noin{\noindent}
\def\H{{I\!\!H}}

\def\R{{\mathbb{R}}}

\def\l{\lambda}

\def\0{\overline 0}

\def\pn{\par\noindent}

\def\psn{\par\smallskip\noindent}

\font\tenmsb=msbm10 \font\sevenmsb=msbm7 \font\fivemsb=msbm5
\newfam\msbfam
\textfont\msbfam=\tenmsb \scriptfont\msbfam=\sevenmsb
\scriptscriptfont\msbfam=\fivemsb
\font\teneufm=eufm10 \font\seveneufm=eufm7 \font\fiveeufm=eufm5
\newfam\eufmfam
\textfont\eufmfam=\teneufm \scriptfont\eufmfam=\seveneufm
\scriptscriptfont\eufmfam=\fiveeufm

\begin{document}

\title{On Minty's theorem in the Heisenberg group}
\author{
 A. Calogero\thanks{
Dipartimento di Matematica e Applicazioni, Universit\`a degli Studi
di Milano--Bicocca, Via Cozzi 53, 20125 Milano, Italy ({\tt
andrea.calogero@unimib.it}, $^+$corresponding author: {\tt
rita.pini@unimib.it})}\; and
 R. Pini$^{*+}$
 }

\maketitle
\begin{abstract}
\noindent In this paper we extend some classical results of Convex
Analysis to the sub-Riemannian setting of the Heisenberg group. In
particular, we provide a horizontal version of Minty's theorem
concerning maximal H-monotone operators defined in the Heisenberg
group with values in the first layer of its Lie algebra.
\end{abstract}

\noindent {\bf Key words}: Heisenberg group, convexity,
subdifferential, monotone map, cyclic monotone map, maximal
monotonicity \vskip0.1truecm \noindent {\bf MSC}: Primary: 26B25
Secondary: 47H05; 53C17

\bigskip

\section{Introduction}

Maximal monotone maps in Euclidean spaces and, more in general, in
Hilbert spaces, play key roles in several settings, specifically in
fixed point theorems and in solving generalized equations. A well
known result, the celebrated Minty theorem, provides a
characterization of maximal monotonicity (see \cite{RoWe2004}):
given a monotone set-valued map $T:X\rightrightarrows X,$ where $X$
is a Hilbert space, then $T$ is maximal monotone if and only if
$I+\lambda T$ is surjective onto $X,$ for every $\lambda
>0;$ in this case, the map $(I+\lambda T)^{-1}$ is single-valued on
$X.$ A remarkable implication of Minty's theorem is the possibility
of approximating $T$ in some sense by single-valued maps
$T_{\lambda},$ called \emph{Yosida approximations}, that are also
maximal monotone and are defined as
$T_{\lambda}=(I-J_{\lambda})/\lambda,$ where $J_{\lambda}=(I+\lambda
T)^{-1}.$

The most notable example of a maximal monotone map arises as the
subdifferential $\partial f$ of a convex function $f:\R^n\to\R$.
More precisely, the subdifferential of the function $f$ at $x$  is
defined as the following subset of $\R^n:$
$$
\partial f(x):=\{v\in \R^n:\; f(y)\ge f(x)+\langle v,
y-x\rangle,\quad \forall y\in \R^n\};
$$
the subdifferential map is therefore given by:
$$
\partial f:\R^n\rightrightarrows \R^n,\qquad x\mapsto \partial f(x).
$$
In this case, if $T=\partial f,$ then the classical Moreau theorem
provides a useful approximation of $f$ via convex and Fr\'echet
differentiable functions $f_{\lambda},$ converging upward to $f$,
and whose gradient is given by the Yosida approximation $(\partial
f)_{\lambda}.$

The aim of this paper is to start the investigation of a possible
extension of these results, that are peculiar of Convex Analysis, to
the sub-Riemannian background, starting from the simplest situation
of the Heisenberg group $(\H,\circ)$. In this framework, a first
step towards this study makes use of the concept of convexity and
subdifferentiability of a function $u:\H\to \R$. In the last few
years several notions of convexity have been introduced in the
Heisenberg group and, more generally, in Carnot groups, but the
notion of horizontal convexity (H-convexity) turned out to be the
most suitable to match their sub-Riemannian structure (see
\cite{DaGaNh2003}, \cite{LuMaStr2004}). The naturally associated
H-subdifferential of a function $u$ at a point $g\in \H$  is the
(possibly empty) subset of the first layer $V_1$ of the Lie algebra
of $\H$:
$$
\partial_Hu(g):=\{p\in V_1:\, u(g\circ\exp w)\ge u(g)+\langle w,p\rangle,
\;\forall w\in V_1\},
$$
where $\circ$ denotes the group law in $\H,$ and $V_1$ can be
identified with $\R^2$ (see Section 2 for the details). The
H-subdifferential map is therefore the set-valued map
$$
\partial_Hu:\H\rightrightarrows V_1,\qquad g\mapsto \partial_Hu(g).
$$
It is worthwhile noticing that, while in the Euclidean case the
subdifferential at a point can be identified with a subset of the
same space $X^*=X,$ in the framework of the Heisenberg group there
is a reduction in dimension passing from $\H$ to $V_1;$ this fact
will give rise to some pathological effects. \par\smallskip The
notion of H-monotonicity, that turns out to be a particular case of
the more general notion of $c$ H-monotonicity (see \cite{CaPi2012}),
fits the monotonicity of maps in Euclidean spaces to the horizontal
structure of $\H.$ As in the Euclidean case, the notion of
H-convexity for functions is strictly related to the notion of
H-monotonicity of the H-subdifferential map. In Section 3 we show
that the map $\partial_Hu,$ if $u:\H\to\R$ is H-convex, is maximal
H-monotone.

More generally, we consider a map $T:\H\rightrightarrows V_1,$ which
is maximal monotone in the horizontal sense and try to see whether a
suitable version of Minty's type theorem still holds. In Section 4,
we prove our main result, that is, a horizontal version of Minty's
theorem in a slightly weaker form:
\begin{theorem}\label{main result}
Let $T:\H\rightrightarrows V_1$ be an H-monotone operator with
$\mathrm{dom}(T)=\H.$ If $T$ is maximal H-cyclically monotone, then
the map $(\xi_1+\lambda T)|_{H_g}$ is surjective onto $V_1$ for
every $g\in \H.$

Conversely, if the map $(\xi_1+\lambda T)|_{H_g}$ is surjective onto
$V_1$ for every $g\in \H,$ then $T$ is maximal H-monotone.
\end{theorem}
Here $H_g$ denotes the horizontal plane associated to the point
$g\in\H,$ i.e. the set of all the points in $\H$ that are reachable
from $g$ via horizontal segments (see Section 2 for the details).
Nevertheless, despite one is lead to foresee a possible parallelism
with the Euclidean case, this one fails when trying to show that,
for a fixed $g\in\H,$ the map $(\xi_1+\lambda
T)^{-1}:V_1\rightrightarrows H_g$ is single-valued. As a matter of
fact, we show (see Example \ref{example2}) that this is no longer
true even in the case $T=\partial_Hu,$ where $u$ is an H-convex
function on $\H.$ This fact is quite unexpected, for, in this case,
there is no change in dimension when we pass from $V_1$ to the
horizontal plane $H_g,$ whose topological dimension is 2.

\section{Preliminaries}
The Heisenberg group $\H$ is the Lie group given by the underlying
manifold $\R^3$ with the non commutative group law
$$
g\circ g'=(x,y,t)\circ(x',y',t')= \left(x+x',y+y',t+t'+2(x'y-xy')\right),
$$
unit element $e=(0,0,0),$ and $g^{-1}=(-x,-y,-t)$. Left translations
and anisotropic dilations are, in this setup, $L_{g_0}(g)=g_0\circ
g$ and $\delta_\lambda(x,y,t)=\left(\lambda x,\lambda
y,\lambda^2t\right).$
\vskip 0.1truecm The  differentiable  structure on $\H$ is
determined  by  the left invariant vector fields
$$
X=\partial_{x}+2y\partial_{t},\qquad
Y=\partial_{y}-2x\partial_{t},\qquad T=\partial_{t},\qquad
\texttt{\rm with}\ \ [X,Y]=-4T.
$$
\noin The vector field $T$ commutes with the vector fields $X$ and
$Y$; $X$ and $Y$ are called \it horizontal vector fields\rm.

The Lie algebra $\mathfrak{h}$ of $\H$ is the stratified algebra
$\mathfrak{h}=\R^3=V_1\oplus V_2,$ where $V_1={\rm
span}\left\{X,Y\right\},$ $V_2={\rm span}\left\{T\right\};$
$\langle\cdot ,\cdot\rangle$ will denote the inner product. Via the
exponential map $\exp:\mathfrak{h}\to\H$ we identify the vector
$\alpha X+\beta Y+\gamma T$ in $\mathfrak{h}$ with the point
$(\alpha, \beta, \gamma)$ in $\H;$ the inverse $\xi:
\H\to\mathfrak{h}$ of the exponential map has the unique
decomposition $\xi=(\xi_1,\xi_2),$ with $\xi_i:\H\to V_i.$ Since we
identify $V_1$ with $\R^2$ when needed, $\xi_1:\H\to V_1\sim\R^2$ is
given by $\xi_1(x,y,t)=(x,y).$ We say that $\gamma:[0,1]\to\H$ is a
\emph{horizontal segment} if
$\gamma(\lambda)=g\circ\delta_\lambda(\exp w)=g\circ\exp(\lambda w)$
for some $g\in\H$ and $w\in V_1$ fixed. For more details on the
structure of the Heisenberg group see, e.g., \cite{CaDaPaTy2007} and
\cite{BoLaUg2007}.
\vskip 0.1truecm The main issue in the analysis on the Heisenberg
group is that the classical differential operators are considered
only in terms of the horizontal fields. For any open subset $\Omega$
of $\H,$ let us denote by $\Gamma^1(\Omega)$ the class of functions
having continuous derivatives with respect to the vector fields $X$
and $Y.$ We recall that the horizontal gradient of a function
$u\in\Gamma^1(\Omega)$ at $g\in\Omega$ is the $2$--vector
$$
(\nabla_H u)(g)= \left((Xu)(g), \ (Yu)(g)\right),
$$
written with respect to the basis $\{X,Y\}$ of $V_1.$

The notion of horizontal subdifferential of a function $u:\H\to \R$
at a point $g$ takes into account the sub-Riemannian structure of
$\H.$ This horizontal structure relies on the notion of horizontal
plane: given a point $g_0\in\H$, the {\sl horizontal plane}
$H_{g_0}$ associated to $g_0$ is the plane in $\H$ defined by
\begin{align*}
H_{g_0}=L_{g_0}\left(\exp(V_1)\right)&=\left\{g=(x,y,t):
t=t_0+2y_0x-2x_0y\right\}\\
&=\left\{g_0\circ\exp w,\, w\in V_1\right\}.
\end{align*}
We note that  $g'\in H_{g}$ if and only if $g\in H_{g'}.$

\begin{definition}\label{Def
Hsubdiffer1} Let $u:\Omega\to\R,$ where $\Omega$ is a subset of
$\H.$ The {\sl horizontal subdifferential} (or H--subdifferential)
of $u$ at $g\in\Omega$ is the set
\begin{equation}\label{Hsub}
\partial_Hu(g)=\{ v\in V_1:\ u(g\circ\exp w)\ge u(g)+\langle v,
w\rangle, \quad \forall w\in V_1:\ g\circ\exp w\in\Omega\}.
\end{equation}
\end{definition}
If $v\in\partial_H u(g),$ we say that $v$ is an {\sl H--subgradient}
of $u$ at $g.$ The set-valued map
$$
\partial_Hu:\H\rightrightarrows V_1,\qquad g\mapsto
\partial_Hu(g)
$$
is called the \emph{H-subdifferential} of $u$.

Let us recall that, given a set-valued map $T:\H\rightrightarrows
V_1,$ its domain is defined as $\mathrm{dom}(T):=\{g\in\H:\,
T(g)\neq \emptyset\}, $ and its graph is the subset of $\H\times
V_1$ given by $ \mathrm{gph}(T):= \{(g,v):\, g\in \mathrm{dom}(T),\,
v\in T(g)\}. $

Convex functions in the Heisenberg group setting were first
introduced by Luis Caffarelli (in unpublished work from 1996). This
notion did not really surface in the literature until 2002, when it
was independently formulated and studied, in the more general
setting of Carnot groups, in \cite{LuMaStr2004} and in
\cite{DaGaNh2003}. Essentially, a convex function in $\H$ is a
function whose restriction to horizontal segments are Euclidean
convex functions of one variable.
\begin{definition}\label{def function H-convex}{\rm
(H-convexity)} A function $u:\Omega\subset \H\to \R$ is called {\it
H-convex} if
$$
u(g\circ\exp(\l v))\le u(g)+\lambda\left(u(g\circ\exp v)-u(g)\right)
$$
for every $g\in\Omega,$ $v\in V_1$ and $\l\in [0,1],$ provided
$g\circ\exp(\l v)\in \Omega$.
\end{definition}
Despite the notion of H-convexity requires a suitable behavior only
on the horizontal segments, H-convex functions enjoy some nice
regularity properties, as local Lipschitz continuity (see
\cite{BaRi2003}) and hence differentiability almost everywhere in
horizontal directions.

Likewise the Euclidean case, the relation between H-convexity and
H-subdifferentiability of a function is not unexpected. The
investigation of this connection started in the pioneering work
\cite{DaGaNh2003}, and it was carried out in \cite{CaPi2011}: we
mention that in the recent paper \cite{BaDr2013} the authors study
this connection without involving any group structure. In our
context, the following result holds:
 \begin{theorem}{\rm (see Proposition 10.5 in
\cite{DaGaNh2003} and Theorem 4.4 in \cite{CaPi2011})} Let $\Omega$
be an open subset of $\H,$ and $u:\Omega\to \R.$ Then $u$ is
H-convex if and only if $\partial_Hu(g)\neq \emptyset,$ for every
$g\in \Omega.$
\end{theorem}
It is worthwhile noticing that if $u$ is H-convex and $u\in
\Gamma^1(\Omega),$ then $\partial_H u(g)=\{\nabla_H u(g)\}$ for
every $g\in \Omega.$

\section{On the H-subdifferential map of an H-convex function}

The main result of this section, Theorem \ref{propo 1}, concerns the
relationship between H-convex functions and the monotonicity
properties of their subdifferential maps.

The notions of H-convexity for a function $u:\H\to\R,$ and of the
H-subdifferential map $\partial_H u:\H\rightrightarrows V_1\sim
\R^2$ given in the previous section, supply the guidelines in order
to introduce and to study the notion of monotonicity for set-valued
maps $T:\H\rightrightarrows V_1.$

First of all, consider an H-convex function $u:\H\to\R.$ For every
$g\in\H$ and $g'\in H_g$ we obtain, by \eqref{Hsub},
\begin{equation}\label{con 1}
u(g)\ge u(g')+
\langle v,\xi_1(g')-\xi_1(g)\rangle,
\end{equation}
for every $v\in \partial_H u(g).$ Since $g\in H_{g'},$ on the other
hands we obtain
\begin{equation}\label{con 2}
u(g')\ge u(g)+
\langle v',\xi_1(g)-\xi_1(g')\rangle,
\end{equation}
for every $v'\in \partial_H u(g').$ Keeping the comparison between
Euclidean notions and horizontal ones, we provide the following
\begin{definition}\label{definizione monotone}{\rm (H-monotonicity)}
Given a set-valued map $T:\H\rightrightarrows V_1,$ $T$ is said to
be H-monotone if
$$
\langle v'-v,\xi_1(g')-\xi_1(g)\rangle\ge 0,
$$
for every $g\in \mathrm{dom}(T),$ $g'\in \mathrm{dom}(T)\cap H_g,$
$v\in T(g),$ $v'\in T(g').$
\end{definition}

In order to introduce the notion of H-cyclic monotonicity, we say
that the set $\{g_i\}_{i=0}^n\subset \H,$ $n>0,$ is an
\emph{H-sequence} if $g_{i+1}\in H_{g_i},$ for every $i=0,\dots,
n-1.$  An H-sequence is \emph{closed} if $g_n\in H_{g_0};$ in this
case, we usually set $g_{n+1}=g_0.$ In particular, if we consider an
H-convex function as before, and a closed H-sequence
$\{g_i\}_{i=0}^n,$ we obtain
\begin{equation}\label{con 3}
u(g_{i+1})\ge u(g_i)+ \langle v_i,\xi_1(g_{i+1})-\xi_1(g_i)\rangle,
\end{equation}
for every $v_i\in \partial_H u(g_i),\ i=1,\ldots, n.$ As in the
Euclidean framework, adding up the two sides of all the previous $n$
inequalities we are lead to the following notion (see Definition 6.1
in \cite{CaPi2012}):
\begin{definition}\label{H-cyclical monotonicity}
{\rm (H-cyclic monotonicity)} We say that $\mathcal{R}\subset
\H\times V_1$ is an H-cyclically monotone set if, for every sequence
$\{(g_i,v_i)\}_{i=0}^{n}\subset \mathcal{R}$ such that
$\{g_i\}_{i=0}^{n}$ is a closed H-sequence, we have that
\begin{equation}\label{def c Hcyclically}
\sum_{i=0}^{n} \langle \xi_1(g_{i+1}),v_i\rangle \le \sum_{i=0}^{n}
\langle \xi_1(g_{i}),v_i\rangle.
\end{equation}
A set-valued map $T:\H\rightrightarrows V_1$ is an H-cyclically
monotone map if $\mathrm{gph}(T)$ is H-cyclically monotone.
\end{definition}
It is clear that, for an H-convex function $u,$ the inequalities
\eqref{con 1} and \eqref{con 2} entail that $\partial_H
u:\H\rightrightarrows V_1$ is H-monotone; moreover, summing up the
$n$ inequalities in \eqref{con 3}, we get that $\partial_H u$ is
H-cyclically monotone. In the following we will provide an example
of an H-monotone map that is not H-cyclically monotone (see Example
\ref{example1}).

A first study of the properties of H-monotonicity appears in
\cite{CaPi2011} where the authors prove (see Theorem 6.4) that if
$u:\H\to\R$ is H-convex, then
\begin{equation}\label{Rock function}
u(g)=u(g_0)+\sup_{\cal P} \left\{\sum_{i=0}^{n-1} \langle
v_i,\xi_1(g_{i+1})-\xi_1(g_i)\rangle\right\},
\end{equation}
where $g_0\in H$ is fixed and
$$ {\cal P}=\Bigl\{\{(g_i,v_i)\}_{i=0}^{n}\subset \mathrm{gph}(\partial_hu),\, \{g_i\}_{i=0}^{n}\ \texttt{\rm closed
H-sequence},\, n>0\Bigr\};
$$
the right hand side of \eqref{Rock function} is called
\emph{Rockafellar function}. In the proof of this result the
H-cyclic monotonicity of $\mathrm{gph}(\partial_H u)$ is crucial.
The role of the horizontal version of the Rockafellar function and
the connection between  H-cyclically monotone sets and H-convex
functions is  emphasized in a subsequent result in
 \cite{CaPi2012}.  A useful version of Theorem 6.6 in
\cite{CaPi2012} in our context is the following:

\begin{theorem}\label{theorem rockafellar c}
Let $T:\H\rightrightarrows V_1$ be an H-cyclically monotone map with
$\mathrm{dom}(T)=\H.$ Then there exists an H-convex function
$u:\H\to \R $ such that
\begin{equation}\label{inclusione teo rock}
\texttt{\rm gph}(T)\subset \texttt{\rm gph}(\partial_H u).
\end{equation}
\end{theorem}
We emphasize that the function $u$ mentioned in the previous result
is a Rockafellar function.

It is well known that in the Euclidean case the notion of maximality
is crucial. We say that
 the set-valued map $T:\H\rightrightarrows V_1$ is \emph{maximal H-monotone} (\emph{maximal
H-cyclically monotone}) if there are no H-monotone (H-cyclically
monotone) set-valued maps $T':\H\rightrightarrows V_1$ such that
$T(g)\subset T'(g)$ for every $g\in \H,$ and $T(g')\subsetneq
T'(g'),$ for some $g'\in \H.$

Notice that any H-cyclically monotone map that is maximal
H-monotone, is maximal H-cyclically monotone, i.e. there are no
H-cyclically monotone maps $T'$ such that $\mathrm{gph}(T)\subsetneq
\mathrm{gph}(T').$

In the following the next remark will be of some use:
\begin{remark}
Let $T:\H\rightrightarrows V_1$ be H-monotone. Then $T$ is maximal
H-monotone if and only if for all $(g,v)\notin \mathrm{gph}(T),$
there exists $g'\in H_g$ and $v'\in T(g')$ such that
$$
\langle v-v',\xi_1(g)-\xi_1(g')\rangle<0.
$$
\end{remark}

The main result of this section is the following

\begin{theorem}\label{propo 1}
\quad

\begin{itemize}
\item[i.] If $u:\H\to \R$ is an H-convex function, then the set-valued map $T=\partial_H u:\H\rightrightarrows V_1$
is maximal H-monotone; since it is H-cyclically monotone, it is also
maximal H-cyclically monotone.

\item[ii.] If $T:\H\rightrightarrows V_1$ is a maximal H-cyclically monotone with $\mathrm{dom}(T)=\H,$ then there exists an H-convex function $u$ such that
$T=\partial_H u$.
\end{itemize}

\end{theorem}

\pn \textbf{Proof:} i. Let $u:\H\to \R$ be H-convex. First of all
notice that, from Theorem 4.4 in \cite{CaPi2011},
$\mathrm{dom}(T)=\H.$ We have just seen at the beginning of this
section that the monotonicity of $\partial_Hu$ follows trivially
from the H-convexity of $u.$ Straightforward computations show that
$\partial_Hu$ is also H-cyclically monotone, according to Definition
\ref{H-cyclical monotonicity}.

Let us prove the maximal monotonicity. Fix $g\in \H,$ and take any
$v\notin \partial_Hu(g);$ by the definition of horizontal
subgradient, there exists $g'\in H_g,$ $g'=g\circ \exp z$ for some
$z\in V_1,$ such that
\begin{equation}\label{absurd}
u(g')<u(g)+\langle v,\xi_1(g')-\xi_1(g)\rangle.
\end{equation}
Set $\phi(\cdot)=u(\cdot)-\langle v,\xi_1(\cdot)-\xi_1(g)\rangle;$
we have that
\begin{equation}\label{dim 1}
\phi(g')<\phi(g).
\end{equation}
From the H-convexity of $\phi:\H\to\R,$ it is clear that the
restriction of $\phi$ to the horizontal segment $[g,g']\subset
H_g\cap H_{g'},$ i.e. the function $\l \to \phi(g\circ\exp (\l z)),$
$\l\in[0,1],$ is Euclidean convex. Let us notice that
$$
\phi(g\circ\exp(\lambda z)\circ \exp(\lambda' z))=\phi(g\circ\exp((\lambda+\lambda') z)).
$$
Then, from \eqref{dim 1}, there exists $\tilde\l \in (0,1)$ such
that, for $\tilde{g}=g\circ\exp(\tilde{\l}z)\in H_g,$
$$
\lim_{\lambda\to 0^+}\frac{\phi(\tilde{g}\circ\exp(\lambda
z))-\phi(\tilde{g})}{\lambda}<0.
$$

Let us denote by $\phi'(\tilde{g};z)$ the previous limit. Now take
any $\tilde{v}\in
\partial_Hu(\tilde{g}).$ Then Proposition 4.1 in \cite{CaPi2011}
implies that
$$
u'(\tilde{g};v)\ge \langle \tilde{v},v\rangle,\qquad \forall v\in
V_1; $$ in particular, $u'(\tilde{g}; z)\ge \langle
\tilde{v},z\rangle.$ Moreover, from the definition of $\phi,$
$$
u'(\tilde{g};z)=\phi'(\tilde{g};z)+\langle v,z\rangle<\langle
v,z\rangle,
$$
and therefore, recalling that $\tilde\l>0,$
$$
\tilde\l\langle \tilde{v}-v, z\rangle=\langle \tilde{v}-v,
\xi_1(\tilde{g})-\xi_1(g)\rangle <0,
$$
contradicting the H-monotonicity of $\partial_H u.$

ii. Suppose now that $T$ is maximal H-cyclically monotone with
$\mathrm{dom}(T)=\H.$ Then, Theorem \ref{theorem rockafellar c}
implies that there exists an H-convex function $u$ such that
$T(g)\subset
\partial_H u(g),$ for every $g\in\H.$ From the maximal H-cyclic monotonicity of $T$ and i., we get that $T=\partial_Hu.$
 \QED

As in \cite{AubFra1990}, Example 12.7, the following result holds in
the Heisenberg setting:
\begin{proposition}
Let $T:\H\rightrightarrows V_1$ be a continuous H-monotone and
single-valued map with $\mathrm{dom}(T)=\H.$ Then $T$ is
 maximal H-monotone.
\end{proposition}
\textbf{Proof:} Let $g\in\H$ and $v\in V_1$ such that
$$
\langle v -T(g'),\xi_1(g)-\xi_1(g')\rangle \ge 0,\qquad \forall g'\in
H_g;
$$
we will prove that $v=T(g).$ Set $g'=g\circ
\exp(-\lambda(\xi_1(\tilde g)-\xi_1(g))),$ for some $\tilde g\in\H$
and $\lambda>0;$ we obtain, from the previous inequality and
dividing by $\lambda,$
$$
\langle v -T(g\circ\exp(-\lambda(\xi_1(\tilde g)-\xi_1(g)))),\xi_1(\tilde g)-\xi_1(g)\rangle \ge
0,\qquad \forall \lambda>0.
$$
Taking the limit as $\lambda\to 0^+,$ we obtain
$$
\langle v -T(g),
\xi_1(\tilde g)-\xi_1(g)\rangle \ge 0;
$$
the generality of $\tilde g$ implies that $v=T(g).$
 \QED

\begin{remark}
Let $\tilde T:\R^2\rightrightarrows\R^2$ and $ T:\H\rightrightarrows
V_1=\R^2$ be maps with $\mathrm{dom}(\tilde T)=\R^2,\
\mathrm{dom}(T)=\H,$ such that
 \begin{equation}\label{maps} T(x)=\tilde
T(\xi_1(x)),\qquad\forall x\in\H.
\end{equation}
If $\tilde T$ is cyclically monotone, then $T$ is H-cyclically
monotone. The converse is false, in general.
\end{remark}

This difference between the \lq\lq Euclidean" cyclic monotonicity of
$T$ and the H-cyclic monotonicity of $\tilde T$ is delicate, as the
next example shows.
\begin{example}\label{example1}
{\rm Let us consider $\tilde T:\R^2\to\R^2$ defined by $\tilde
T(x,y)=Q(x,y)^T,$ where $Q$ is a $2\times 2$ matrix. It is
well-known (see \cite{Ro1969}, p. 240) that $\tilde T$ is monotone
if and only if $\frac{1}{2}(Q+Q^T)$ is positive semidefinite;
moreover $\tilde T$ is cyclically monotone if and only if $Q$ is
symmetric and positive semidefinite.

Let us consider the particular case
$Q=\left(\begin{array}{cc}3&2\\-2&4\end{array}\right):$ in this
case, $\tilde T$ is maximal monotone (note that $\tilde T$ is
continuous). Define the map $T:\H\rightrightarrows V_1=\R^2$ as in
\eqref{maps}: explicitly,
$$
T(x,y,t)=\tilde T (x,y)=(3x+2y,-2x+4y).
$$
Consider now the function $u:\H\to\R$ defined by
$$
u(x,y,t)=\frac{3}{2}x^2+2y^2+t.
$$
It is clear that $u$ is Euclidean convex, hence H-convex; its
regularity implies that $\partial_H u(g)=\{\nabla_H u(g)\}=T(g),$
for every $g=(x,y,t)\in \H.$ By Theorem \ref{propo 1}, $T$ is
maximal H-cyclically monotone.

On the other hand, if we consider the matrix
$Q=\left(\begin{array}{cc}3&0\\-2&4\end{array}\right),$ the related
operator $\tilde T$  is again maximal monotone. The map $ T$ is
defined by
$$
T(x,y,t)=\tilde T (x,y)=(3x,-2x+4y).
$$
Let us prove that $T$ is maximal H-monotone, but not maximal
H-cyclically monotone. The maximality follows from the continuity of
$T;$ besides, the H-monotonicity can be inferred from the
monotonicity of $\tilde T.$ Suppose that $T$ is H-cyclically
monotone: then, $\mathrm{dom}(T)=\H$ and Theorem
 \ref{theorem rockafellar c} imply that there exists an H-convex function
$u:\H\to\R$ such that $\texttt{\rm gph}(T)\subset\texttt{\rm
gph}(\partial_H u);$ the maximality of $T$ gives that $T=\partial_H
u.$ Since $\partial_H u(g)$ is a singleton for all $g\in\H,$ then
$u$ is Pansu differentiable within $\H$ and $\partial_H
u(g)=\{\nabla_H u(g)\}$ for every $g\in \H$ (see Theorem 1.3 in
\cite{MaSci2010}). Hence, we have that
$$
\nabla_H u(g)=(u_{x}+2yu_{t},u_{y}-2x
u_{t})=(3x,-2x+4y)=T(g).
$$
An easy computation shows that a function $u$ satisfying the
previous equality does not exist, hence  $T$ is not  H-cyclically
monotone.}
\end{example}

\section{The main result}
This section is devoted to the more subtle result, where we prove
the horizontal version of the Minty theorem contained in Theorem
\ref{main result}.

In order to prove the assertion, we need the following two lemmata:

\begin{lemma}\label{superlinearity}
Let $u:\H\to\R$ be an H-convex function, and $g_0\in \H.$ Consider
the function $f:H_{g_0}\to\R$ defined by
$$
f(g)=u(g)+\frac{|\xi_1(g)|^2}{2}.
$$
Then, $\displaystyle \lim_{|\xi_1(g)|\to \infty}f(g)=+\infty. $
\end{lemma}

\pn \textbf{Proof:} Take $v\in \partial_Hu(g_0);$ then, for any
$g\in H_{g_0},$
$$
u(g)\ge u(g_0)+\langle v,\xi_1(g)-\xi_1(g_0)\rangle=
 \langle v,\xi_1(g)\rangle +u(g_0)-\langle v,\xi_1(g_0)\rangle.
$$
Hence, for $g\in H_{g_0}$ and $c=u(g_0)-\langle
v,\xi_1(g_0)\rangle,$
$$
f(g)\ge \langle v, \xi_1(g)\rangle +c +\frac{|\xi_1(g)|^2}{2},
$$
therefore we get the assertion. \QED

\begin{lemma}\label{left traslation}
Let $T:\H\rightrightarrows V_1.$ Let us consider a fixed point
$g_0\in\H$ and the map $T_0(g)=T(g_0\circ g).$ Then,
\begin{itemize}
\item[i.] $T$ is H-monotone if and only if $T_0$ is H-monotone;

\item[ii.]
$ \mathrm{rge}(T_0+\xi_1)|_{H_g}=V_1,$ $\forall g\in \H, $ if and
only if $\mathrm{rge}(T+\xi_1)|_{H_g}=V_1,$ $\forall g\in \H.$
\end{itemize}

\end{lemma}

\pn \textbf{Proof:} i. Let $T$ be  H-monotone. Suppose that $v\in
T_0(g),$ $g'\in H_g,$ $v'\in T_0(g');$ then we have that $v\in
T(g_0\circ g),$ $g_0\circ g'\in H_{g_0\circ g}$ and $v'\in
T(g_0\circ g).$ In addition, $\xi_1(g)-\xi_1(g')=\xi_1(g_0\circ
g)-\xi_1(g_0\circ g').$ Therefore,
$$
\langle v-v',\xi_1(g)-\xi_1(g')\rangle \ge 0,
$$
from the H-monotonicity of $T.$ Hence $T_0$ is H-monotone. The
converse can be analogously proved. \psn ii. Let us remark that
\begin{eqnarray*}
\mathrm{rge}(T+\xi_1)|_{H_g}&=&\mathrm{rge}\bigl(T(g_0\circ (g_0^{-1}\circ\cdot))+\xi_1(\cdot)\bigr)|_{H_g}\\
&=&\mathrm{rge}(T(g_0\circ(\cdot))+\xi_1(\cdot)+\xi_1(g_0))|_{H_{g_0^{-1}\circ g}}\\
&=&\mathrm{rge}(T_0+\xi_1)|_{H_{g_0^{-1}\circ g}}+\xi_1(g_0).
\end{eqnarray*}
\QED The proof of Theorem \ref{main result} is based on the
following comparison lemma for the horizontal normal map (see
Theorem 3.1 in \cite{BaCaKr2013}), that is quoted below in a
simplified version:
\begin{theorem}\label{comparison} Let $\Omega_0\subset
\H$ be an open and bounded set, and $u,v:\H\to\R$ be $H$-convex
functions. Let $g_0\in \Omega_0$ such that $u(g_0)\le v(g_0),$ and
$u\ge v$ on $\partial \Omega_0\cap H_{g_0}.$ If $v\in
\partial_Hv(g_0)$ satisfies the inequality
$$
v(g)>v(g_0)+\langle v,\xi_1(g)-\xi_1(g_0)\rangle,\qquad \forall g\in
\partial \Omega_0\cap H_{g_0},
$$
then $v\in \partial_Hu(\Omega_0\cap H_{g_0}).$
\end{theorem}

We are now in a position to prove our main result.
 \psn  \textbf{Proof
of Theorem \ref{main result}:} Let $T:\H\rightrightarrows V_1$ be a
maximal H-cyclically monotone map with $\mathrm{dom}(T)=\H.$ Theorem
\ref{propo 1} $ii.$ gives that
 there exists a convex function $u:\H\to \R$
such that $T=\partial_H u.$

First of all notice that it is enough to show that, for every
$g_0\in \H,$ there exists $g\in H_{g_0}$ such that $0\in
\partial_Hu(g)+\xi_1(g).$
Indeed, suppose that this assertion is proved. Take any $v\in V_1$
and consider the function
$$
u_v(g)=u(g)-\langle v,\xi_1(g)\rangle.
$$
The function $u_v$ is still H-convex, then there exists $g\in
H_{g_0}$ such that ${0}\in
\partial_Hu_v(g)+\xi_1(g).$
Since $\partial_H(-\langle v,\xi_1(\cdot)\rangle)=\{-v\},$ for every
$g\in \H,$ we have the equality
$$
\partial_H u_v=\partial_H u-v;
$$
therefore, ${0}\in \partial_Hu(g)-v+\xi_1(g),$ or, equivalently,
$v\in
\partial_Hu(g)+\xi_1(g).$

Let us fix $g_0\in\H$ and consider the H-convex function $\phi:\H\to
\R$ defined by
$$
\phi(g)=u(g)+\frac{|\xi_1(g)|^2}{2}.
$$
By Lemma \ref{superlinearity}, $\phi(g)\to +\infty$ whenever
$|\xi_1(g)|\to \infty$ with $g\in H_{g_0}.$ This implies that, for
every $M>0,$ there exists $R>0$ such that
\begin{equation}\label{ineq_on_boundary}
\phi(g)\ge \phi({g_0})+M,\qquad \forall g\in \partial B({g_0},R)\cap H_{g_0},
\end{equation}
where $B(g,r)$ is the open Euclidean ball of center $g\in\H$ and
radius $r.$ Let us denote by $\tau:=\phi({g_0})+M$ and by
$V:\R^{3}\to\R$ the classical Euclidean convex function whose graph
in $\H\times \R$ is the upside-down cone such that the vertex is the
point $({g_0},\phi({g_0}))\in\H\times \R,$ and the basis is
$\partial B({g_0},R)\times\tau\subset  \H\times \R.$ Since  $V$ is
Euclidean convex, clearly it is also H-convex. In order to apply
Theorem \ref{comparison}, we observe that
$$
V({g_0})=\phi({g_0})\quad {\rm  and}\quad \phi(g)\geq\tau= V(g)\ {\rm for\ every}\
g\in \partial B({g_0},R)\cap H_{g_0}.
$$
Moreover, we have that $p_0:=0\in\partial_H V({g_0})$ and
$$
V(g)>V({g_0})+\langle p_0,\xi_1(g)-\xi_1({g_0})\rangle, \ {\rm for\ every}\
g\in \partial B({g_0},R)\cap H_{g_0}.
$$
Theorem \ref{comparison} implies that $0\in\partial_H \phi(
B({g_0},R)\cap H_{g_0})$ and we conclude the first part of the
proof.

Now, let $T:\H\rightrightarrows V_1$ be a set-valued H-monotone map,
with domain $\H,$ such that, for every $g_0\in \H,$
$$
\mathrm{rge}(T+\xi_1)|_{H_{g_0}}=V_1.
$$
 We argue by contradiction and we suppose that there
exists $g_0\in \H,$ and $w\notin T(g_0)$ such that, for every $g\in
H_{g_0},$ and $v\in T(g),$
\begin{equation}\label{inequality}
\langle w -v,\xi_1(g_0)-\xi_1(g)\rangle \ge 0. \end{equation} Let
us first consider the case $g_0=0.$ From the assumptions
$\mathrm{rge}(T+\xi_1)|_{H_0}=V_1, $ therefore
\begin{equation}\label{equality}
w=\tilde v+\xi_1(\tilde g),
\end{equation}
for some $\tilde g\in H_0$ and $\tilde v\in T(\tilde g).$ From
 \eqref{equality}, taking $g=\tilde g$ in \eqref{inequality}, we
obtain
$$
-\langle \xi_1(\tilde g),\xi_1(\tilde g)\rangle \ge 0,
$$
i.e., $\xi_1(\tilde g)=0.$ Since $\tilde g\in H_0,$ we deduce that
$\tilde g=0,$ and $w=\tilde v\in T(0),$ contradicting our assumption
on $w.$

Suppose now that $w\notin T(g_0),$ for some $g_0\in \H,$ but
\eqref{inequality} is satisfied for every $g\in H_{g_0},$ $v\in
T(g).$ Let us consider the set-valued map $T_{0}(\cdot)=T(g_0\circ
\cdot)$ as in Lemma \ref{left traslation}; clearly, $w\notin
T_0(0).$ Moreover, for every $g'={g_0}^{-1}\circ g\in H_0$ and $v\in
T_0(g'),$ inequality \eqref{inequality} gives
$$
\langle w -v,\xi_1(0)-\xi_1(g')\rangle = \langle w
-v,\xi_1({g_0})-\xi_1(g)\rangle \ge 0.
$$
Since, from  Lemma \ref{left traslation}, $T_0$ is monotone, from
the first part of the proof we argue that we get a contradiction,
i.e. $w\in T_0(0),$ or, equivalently, $w\in T(g_0).$ \QED

It is an open question for the authors whether, in Theorem \ref{main
result}, the assumption of maximal H-cyclic monotonicity of $T$ can
be weakened to maximal H-monotonicity, as in the classical case.

\subsection{The open question of the resolvent $J_{\lambda}.$}
In the setting of Minty's theorem, given a maximal monotone
set-valued map $T:X\rightrightarrows X$ where $X$ is a Hilbert
space, for any $\lambda
>0$ one can define a map $J_{\lambda},$ called \emph{resolvent} of
$T$ and given by
$$
J_{\lambda}:X\rightrightarrows X,\quad J_{\lambda}:=(I+\lambda
T)^{-1}.
$$
This map turns out to be non-expansive and single-valued from $X$ to
$X$ (see, for example, \cite{RoWe2004}); furthermore, the
 \emph{Yosida approximation} $T_{\lambda}=(I-J_{\lambda})/\lambda$ is maximal
 monotone as well.
If we consider, in particular, the map $T=\partial \varphi,$
 where $\varphi:X\to [-\infty,\infty]$ is a proper, convex and lower semicontinuous function, a well-known result due to Moreau states that the function
$$
\varphi_\lambda(x)=\min_{x'\in
H}\left(\frac{1}{2\alpha}\left|x'-x\right|+\varphi(x')\right),\quad
\lambda >0
$$
is convex and Fr\'echet differentiable, with gradient
$\varphi_\lambda'=T_\lambda.$ In addition, $\varphi_\lambda(x)$
converges to $\varphi(x)$ as $\lambda\downarrow 0$ for each $x\in
H.$

One may wonder whether a similar property is still true in the
Heisenberg setting, by considering the candidate most likeable to
play the role of $J_{\lambda},$ i.e.,
$$
J^H_{\lambda}:=(\xi_1+\lambda T)^{-1}:V_1\rightrightarrows \H.
$$
If we assume that $T$ is maximal H-cyclically monotone, by Theorem
\ref{main result} $J^H_{\lambda}$ is defined on the whole $V_1.$ In
addition, since rge$(\xi_1+\lambda T)|_{H_g}=V_1,$ for every fixed
$g\in \H,$ one may wonder whether
$$
J^H_{\lambda}:V_1\rightrightarrows H_g
$$
is a single-valued map. Unfortunately this property is not
inherited, even if we restrict our attention to the special but
still exhaustive case of the map $T=\partial_Hu,$ where $u$ is a
real-valued H-convex function defined on the whole $\H.$ The main
reason relies on the following fact: once we restrict the map
$J^H_{\lambda}$ to a fixed horizontal plane $H_{g_0}$, all the
properties of the horizontal subdifferential map $g\mapsto
\partial_H u(g),$ $g\in H_{g_0},$ are lost, since, by definition, they concern the
behaviour of the map only at points that can be related two by two
via the condition $g'\in H_g$ or, equivalently, $g\in H_{g'}$ (it
suffices to look at the notion of H-monotonicity, for instance).

Let us provide an example.
\begin{example}\label{example2}
\rm Let us consider the gauge function $N:\H\to \R$ defined as
$$
N(x,y,t)=((x^2+y^2)^2+t^2)^{1/4}.
$$
It is known that this function is H-convex, but it is not Euclidean
convex (see \cite{DaGaNh2003}). The associated horizontal
subgradient map, for every $g =(x,y,z),$ is given by
$$
\partial_HN(g)=
     \begin{cases}
    \overline{B(0,1)} & (x,y,t)=(0,0,0) \\
    \frac{1}{N^3(g)}\left(
    x(x^2+y^2)+yt,y(x^2+y^2)-xt\right) & (x,y,t)\neq (0,0,0).
    \end{cases}
    $$
We will show that there exists $g''\in \H,$ and $g, g'\in H_{g''},$
$g\neq g',$ such that
$$
(\xi_1(g)+\partial_H u(g))\cap
(\xi_1(g')+\partial_H u(g'))\neq \emptyset.
$$
 Set $g'=e=(0,0,0),$ and notice
that, for every $g\in \H, $ with $\xi_1(g)\neq (0,0),$ there exists
$g''\in H_e$ such that $g\in H_{g''}.$ We will prove that
\begin{equation}\label{ex2}
(\xi_1(g)+\lambda\partial_HN(g))\in \lambda
\overline{B(0,1)}=\xi_1(e)+\lambda\partial_HN(e).
\end{equation}
 Straightforward
computations lead to the following ($g=(x,y,t)$):
$$
|\xi_1(g)+\lambda\partial_HN(g)|^2=
(x^2+y^2)\left(1+\frac{\lambda^2}{N^2(g)}+\frac{2\lambda}{N^3(g)}(x^2+y^2)\right).
$$
Since $g$ can be arbitrarily chosen, we can take, for instance,
$t=1,$ and $x^2+y^2\le 1;$ then, $N(x,y,t)\ge 1,$ and
$$
(x^2+y^2)\left(1+\frac{\lambda^2}{N^2(g)}+\frac{2\lambda}{N^3(g)}(x^2+y^2)\right)\le
(x^2+y^2)(1+\lambda)^2.
$$
Then, if we choose $x,y$ such that $x^2+y^2\le
\lambda^2/(1+\lambda)^2,$ we get the assertion \eqref{ex2}. \psn
\end{example}

As a by-product, the example above brings to light the following
fact:
\begin{remark}{\rm
 There exist H-convex functions $u:\H\to\R$ such that
$$
\partial_Hu(g)\cap \partial_Hu(g')\neq \emptyset,
$$
for some $g,g'$ in the same horizontal plane $H_{g''}.$ }
\end{remark}

\bibliography{capi}
\bibliographystyle{plain}
\end{document}